\documentclass[11pt,a4paper]{article}

\usepackage[T1]{fontenc}
\usepackage[utf8]{inputenc}
\usepackage{lmodern}
\usepackage{microtype}

\usepackage{geometry}
\usepackage{amsmath,amssymb}
\usepackage{booktabs}
\usepackage{graphicx}
\usepackage{float}
\usepackage{xcolor}
\usepackage{url}
\usepackage{hyperref}
\usepackage{siunitx}
\usepackage{pgfplots}
\usepackage{pgfplotstable}
\usepgfplotslibrary{groupplots}
\pgfplotsset{compat=1.18}
\usepackage{lineno}
\usepackage{tabularx}
\usepackage{array}
\usepackage[
backend=biber,
style=numeric,
sorting=nty
]{biblatex}

\newcolumntype{Y}{>{\raggedright\arraybackslash}X}

\hypersetup{
	colorlinks=true,
	linkcolor=black,
	citecolor=blue!60!black,
	urlcolor=blue!60!black,
	pdftitle={Rule-Induced Behavior of Fuzzy Scalar Objective Functions for Reliable Multi-Criteria Decision Making},
	pdfauthor={Olaf Frommann},
	pdfsubject={Rule-Induced Behavior of Fuzzy Scalar Objective Functions for Reliable Multi-Criteria Decision Making},
	pdfkeywords={multi-criteria optimization, scalarization, objective functions, desirability functions, fuzzy logic}
}

\title{Rule-Induced Behavior of Fuzzy Scalar Objective Functions for Reliable Multi-Criteria Decision Making}

\author{
	Olaf Frommann\\
	Institute for Aerospace Technology\\
	Hochschule Bremen, University of Applied Sciences\\
	Bremen, Germany\\
	Corresponding author: \texttt{olaf.frommann@hs-bremen.de}\\
	ORCID: \href{https://orcid.org/0009-0004-2614-4392}{0009-0004-2614-4392}
}

\date{}

\begin{document}
	
	\maketitle
	
	
	
	\begin{abstract}
		Fuzzy-rule-based scalar objective functions provide a flexible way to encode
		qualitative preferences, reference regions, and interactions between criteria in
		multi-criteria optimization and decision making. However, the scalar preference
		landscape induced by such rules can differ substantially from the intended
		decision semantics. This paper investigates how membership placement, implicit
		single-criterion baseline rules, and explicit rule consequents affect the
		behavior of fuzzy scalar objective functions. Two analytically controlled
		bi-criteria Pareto fronts and an embedded two-dimensional dominated-reference
		formulation are used to separate front-selection mechanisms from reference
		improvement behavior. The study shows that apparent reference following can be
		caused by flat rule-activation plateaus, whereas genuine reference following
		requires localized minima with low tie ambiguity. In the dominated-reference
		setting, global memberships with competing consequents recover robust Pareto
		tradeoffs but do not necessarily improve each reference design. By contrast, 
		a reference-based three-class rule set consistently improves
		dominated references, recovers the Pareto set, and avoids plateau-driven
		selection in the present tests. The results provide
		diagnostic metrics and practical guidance for constructing fuzzy scalarizations
		whose optimization behavior is consistent with the intended decision semantics.
	\end{abstract}
	
	
	\section{Introduction}
	\label{sec:introduction}
	
	Many engineering optimization problems require a single design decision in the
	presence of competing criteria. Pareto-front information is useful because it
	describes non-dominated compromises~\cite{Miettinen1999,Ehrgott2005}, but it
	does not by itself define which compromise should be selected. Moreover, the
	Pareto front is often unavailable or impractical to derive due to significant
	computational effort.
	
	Scalarization and preference-based formulations are therefore commonly used to
	convert multi-criteria problems into single-objective optimization
	problems~\cite{MarlerArora2004,Miettinen1999,Steuer1986}. Related scalar
	preference constructions include desirability
	functions~\cite{Harrington1965,DerringerSuich1980}, which also map several
	criteria or responses to a single scalar preference value.
	
	Fuzzy goal programming and fuzzy multiple-objective decision-making methods use
	membership functions to express aspiration levels and satisfaction degrees in
	multi-criteria problems~\cite{BellmanZadeh1970,Zimmermann1978,Sakawa1993,Hannan1981,Tiwari1987,LaiHwang1994}. 
	In many
	applications, these memberships are then aggregated into a scalar satisfaction
	or preference value. Such scalarizations are attractive because they provide a
	direct optimization objective, but they can also hide unintended decision
	mechanisms. In particular, different combinations of membership placement,
	baseline contributions, and rule consequents may produce the same apparent
	selected point for different reasons.
	
	The present study focuses on this induced decision behavior rather than on
	proposing a new fuzzy inference system or a new numerical optimizer. The central
	question is whether a fuzzy scalar objective creates a localized preference
	minimum, a flat plateau, a reference lock, a Pareto tradeoff, or a genuine
	reference-improvement mechanism. This diagnostic perspective complements
	existing fuzzy optimization approaches by analyzing when a fuzzy scalarization
	acts consistently with the intended decision semantics.
	
	The purpose of this paper is to derive diagnostic insight and practical guidance
	for designing fuzzy scalar objective functions for reliable multi-criteria
	decision making.
	Here, reliability means that the scalar objective selects points for the reason
	intended by the rule semantics, rather than because of plateaus, tie-breaking, or
	unintended global attractors. The focus is not on fuzzy logic in general, nor 
	on the performance of a particular numerical optimizer. Instead, analytically 
	controlled test cases are used to isolate the effects of membership structure, 
	implicit single-criterion rules, and explicit rule consequents.
	
	The analysis proceeds in two stages. First, two known Pareto fronts are used to
	study how different fuzzy configurations populate Pareto-front regions. This 
	front-only study reveals mechanisms such as balanced-point collapse,
	genuine reference following, and apparent reference following caused by flat
	rule-activation plateaus. Second, the same Pareto fronts are embedded into a
	two-dimensional design space containing dominated solutions. This dominated-
	reference study tests whether a fuzzy objective can improve a reference design
	rather than merely lock onto it.
	
	The contributions of the paper are:
	\begin{itemize}
		\item a controlled ablation study of fuzzy scalarizations with global and
		reference-based membership placement;
		\item an analysis of how implicit baseline rules and explicit-rules-only
		inference affect the induced scalar preference landscape;
		\item the identification of rule-activation plateaus as a mechanism behind
		apparent reference following and reference stagnation;
		\item a metric-based classification of fuzzy decision behavior, including
		reference stagnation, reference lock, weak reference improvement, Pareto
		reference improvement, and global Pareto tradeoff;
		\item practical guidance for matching fuzzy rule semantics to the intended
		role of the reference point.
	\end{itemize}
	
	\section{Fuzzy scalar objective formulation}
	\label{sec:fuzzy_formulation}
	
	Fuzzy sets were introduced by \cite{Zadeh1965} and were later used for
	decision-making and optimization in fuzzy
	environments~\cite{BellmanZadeh1970,Zimmermann1978,Sakawa1993}. The specific
	scalar objective formulation used in this study is summarized below.
	
	All criteria are treated as normalized minimization criteria,	
	\begin{equation}
		\mathbf{c}=(c_1,\ldots,c_n),
		\qquad c_i\in[0,1],
	\end{equation}
	where smaller values are preferable. For each criterion, fuzzy membership
	values are evaluated for the classes desirable, tolerable, and undesirable:
	$m_{d,i}, \quad m_{t,i}, \quad m_{u,i}.$
	The fuzzy inference produces rule results for the output classes desirable,
	tolerable, and undesirable. Let $r_{d,i},\quad r_{t,i},\quad r_{u,i}$
	denote the activation strengths of the individual desirable, tolerable, and
	undesirable rule results, respectively. The corresponding numbers of rule
	results are denoted by $n_d$, $n_t$, and $n_u$. For compactness, the accumulated
	class strengths are written as
	
	\begin{equation}
		R_d=\sum_{i=1}^{n_d} r_{d,i},
		\qquad
		R_t=\sum_{i=1}^{n_t} r_{t,i},
		\qquad
		R_u=\sum_{i=1}^{n_u} r_{u,i}.
	\end{equation}
	The output evaluation sets used here are triangular sets on the unit output
	axis. With the chosen placement of the desirable, tolerable, and undesirable
	sets, their area centroids are
	\begin{equation}
		s_d=\frac{1}{6},
		\qquad
		s_t=\frac{1}{2},
		\qquad
		s_u=\frac{5}{6}.
	\end{equation}
	The unscaled centroid-based fuzzy value is computed as the weighted mean
	
	\begin{equation}
		F_{\mathrm c}
		=
		\frac{
			s_d R_d + s_t R_t + s_u R_u
		}{
			R_d+R_t+R_u
		}.
		\label{eq:fuzzy_centroid_value}
	\end{equation}
	
	Since $F_{\mathrm c}$ lies in the interval $[1/6,5/6]$ whenever the denominator
	is positive, it is linearly scaled to the unit interval. The minimized fuzzy
	objective is therefore
	
	\begin{equation}
		f_{\mathrm F}
		=
		\frac{F_{\mathrm c}-s_d}{s_u-s_d}
		=
		\frac{3}{2}
		\frac{
			s_d R_d + s_t R_t + s_u R_u
		}{
			R_d+R_t+R_u
		}
		-
		\frac{1}{4}.
		\label{eq:fuzzy_objective}
	\end{equation}
	
	After replacing the output evaluation sets by their fixed centroid scores, the
	rule aggregation is a normalized weighted average of constant consequents. It is
	therefore closely related to zero-order Sugeno or Takagi--Sugeno--Kang fuzzy
	inference, where rule consequents are constants and the final output is obtained
	by normalized rule-strength weighting~\cite{TakagiSugeno1985,SugenoKang1988}.
	Lower values correspond to more preferred outcomes. If the denominator in
	Eq.~\eqref{eq:fuzzy_objective} is zero, the fuzzy objective is undefined for
	that candidate point.
	
	Equivalently, after scaling, the three output classes correspond to the scores
	
	\begin{equation}
		q_d=0,
		\qquad
		q_t=\frac{1}{2},
		\qquad
		q_u=1.
	\end{equation}
	
	Thus, Eq.~\eqref{eq:fuzzy_objective} can also be interpreted as a normalized
	weighted mean of the scaled class scores. The centroid formulation is used 
	here to make the mapping from fuzzy output classes to scalar scores explicit. 
	A detailed description is available in the technical report by
	\cite{Frommann2026OFMCO}.
	
	Two types of membership placement are considered. In global membership
	configurations, desirable and undesirable memberships are defined over the full
	criterion interval. In reference-based configurations, desirable and
	undesirable memberships are defined relative to a reference value. In all
	configurations considered here, tolerable memberships are reference-based.
	
	The present study uses parabolic memberships throughout. The use of only one
	membership type is intentional: the objective is to isolate the effects of
	membership placement, implicit rules, and rule consequents rather than to
	combine these effects with a full comparison of linear, parabolic, and Gaussian
	membership families. All calculations use the publicly available library \texttt{FuzzyGoal}
	in Version~1.2.0~\cite{Frommann2026FuzzyGoal}. Additionally, the programs 
	creating the data for this study are made available for reproducibility 
	purposes~\cite{Frommann2026FuzzyRuleBehavior}.
	
	\subsection{Implicit and explicit rules}
	\label{sec:implicit_explicit_rules}
	
	The fuzzy inference can include implicit single-criterion baseline rules. If
	implicit rules are enabled, each criterion contributes its membership values
	directly to the corresponding output classes. For a bi-criteria problem, the
	implicit contributions have the form
	
	\begin{equation}
		m_{d,1}\rightarrow d,\quad
		m_{t,1}\rightarrow t,\quad
		m_{u,1}\rightarrow u,
		\qquad
		m_{d,2}\rightarrow d,\quad
		m_{t,2}\rightarrow t,\quad
		m_{u,2}\rightarrow u.
	\end{equation}
	
	Thus, even without an explicitly defined interaction rule, each criterion
	provides baseline evidence for the output classes desirable, tolerable, and
	undesirable. In addition, explicitly defined rules may contribute further rule
	results.
	
	If implicit rules are disabled, these baseline contributions are omitted.
	The accumulated output strengths then result only from explicitly defined
	rules. In that case, a candidate point is undefined if no explicit rule produces
	a positive rule result.
	
	This distinction is central to the present study. Implicit rules provide
	continuous single-criterion guidance and often prevent completely flat objective
	regions. However, they can also dominate or mask the effect of explicit
	interaction rules. Explicit-rules-only inference gives the rule base full control over the
	objective, but it can create flat rule-activation plateaus if the active rules
	contribute to only one output class. Such plateaus can then lead to apparent
	reference following through tie-breaking.
	In other words, implicit rules make each criterion individually visible to the
	output aggregation, whereas explicit-rules-only inference evaluates only the
	specified multi-criterion rule base.
	
	\subsection{Rule-activation plateaus}
	\label{sec:rule_activation_plateaus}
	
	A central mechanism in this study is the formation of rule-activation plateaus.
	Such plateaus can occur when the rule base activates only one output class over
	a finite region of the candidate set. In this case, the absolute magnitude of
	the rule activation cancels in the normalized aggregation formula.
	
	To see this, consider a candidate point for which only one output class
	$k\in\{d,t,u\}$ receives a positive accumulated strength,
	\begin{equation}
		R_k=A(\mathbf{c})>0,
		\qquad
		R_l=0
		\quad
		\text{for } l\neq k .
	\end{equation}
	The unscaled centroid value in Eq.~\eqref{eq:fuzzy_centroid_value} then reduces
	to
	\begin{equation}
		F_{\mathrm c}
		=
		\frac{s_k A(\mathbf{c})}{A(\mathbf{c})}
		=
		s_k .
		\label{eq:single_consequent_centroid_plateau}
	\end{equation}
	After the linear scaling in Eq.~\eqref{eq:fuzzy_objective}, the corresponding
	objective value is
	\begin{equation}
		f_{\mathrm F}
		=
		\frac{s_k-s_d}{s_u-s_d}.
		\label{eq:single_consequent_scaled_plateau}
	\end{equation}
	
	Thus, all candidate points for which the same single consequent class is active
	receive exactly the same objective value, independent of the local activation
	strength. The fuzzy objective is flat on that region. If no other rule or
	implicit baseline contribution creates a competing output class, the
	minimization cannot distinguish points within the plateau.
	
	This mechanism is particularly relevant for explicit-rules-only inference with
	single-consequent rule sets. For example, if only the tolerable consequent is
	activated, then
	\begin{equation}
		R_t>0,\qquad R_d=R_u=0,
		\qquad
		F_{\mathrm c}=s_t=\frac{1}{2},
		\qquad
		f_{\mathrm F}=\frac{1}{2}.
	\end{equation}
	
	If only the desirable consequent is activated, then
	$F_{\mathrm c}=s_d=1/6$ and $f_{\mathrm F}=0$. If only the undesirable
	consequent is activated, then $F_{\mathrm c}=s_u=5/6$ and
	$f_{\mathrm F}=1$. In all three cases, the objective value is constant over the
	valid activation region. Points outside the activation region may be undefined
	if no rule contributes positive class strength.
	
	Implicit single-criterion rules can suppress this degeneracy because they add
	baseline contributions to the output classes. Similarly, explicit rule bases
	with competing consequents can break the plateau because the objective then
	depends on the relative strengths of different output classes,
	
	\begin{equation}
		f_{\mathrm F}
		=
		\frac{3}{2}
		\frac{
			s_d R_d + s_t R_t + s_u R_u
		}{
			R_d+R_t+R_u
		}
		-
		\frac{1}{4},
	\end{equation}
	
	rather than on a single class centroid alone. The tie fraction introduced in
	Sec.~\ref{sec:metrics} is therefore used as a diagnostic for detecting whether a
	reported optimum is localized by the fuzzy preference model or selected from a
	flat plateau by tie-breaking.
	
	\subsection{Membership and rule configurations}
	\label{sec:configs_rules}
	
\begin{figure}[ht]
		\centering
		\includegraphics[width=0.7\textwidth]{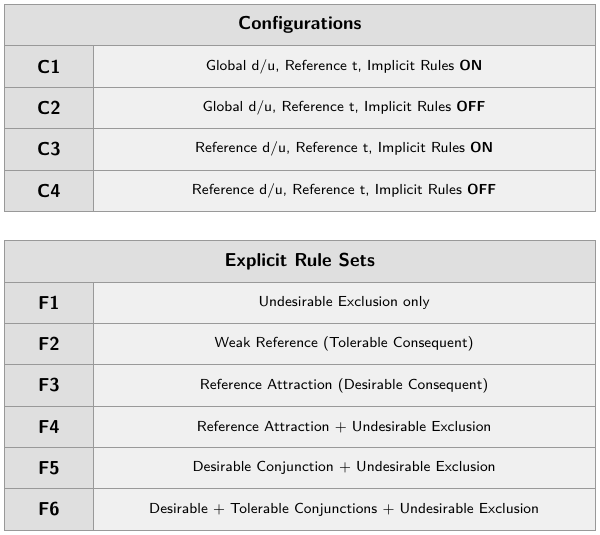}
		\caption{Membership configurations and rule sets used in the ablation study.
			Configurations C1--C4 vary membership placement and implicit-rule policy; rule
			sets F1--F6 vary consequent classes and rule interactions.}
		\label{fig:config_rules}
	\end{figure}

	The four configurations combine membership placement and implicit-rule policy:
	C1 and C2 use global desirable/undesirable memberships, whereas C3 and C4 use
	reference-based desirable/undesirable memberships; C1 and C3 use implicit rules,
	whereas C2 and C4 use explicit rules only. The six explicit rule sets F1--F6
	are summarized in Fig.~\ref{fig:config_rules}.
	
	Rule sets F1--F3 are primarily diagnostic. They isolate exclusion-only,
	weak-reference, and reference-attraction mechanisms. Rule sets F4--F6 are more
	relevant as practical decision templates because they contain competing output
	classes.
	
	\section{Analytical test problems}
	\label{sec:test_problems}
	
	Two analytically defined Pareto fronts are used. They are chosen to distinguish
	supported and non-supported compromise regions while retaining full analytical
	control over the front geometry.
	
	The convex front is
	\begin{equation}
		c_2=g_{\mathrm{conv}}(c_1)
		=
		\frac{4}{3}
		\left[
		\frac{1}{(1+c_1)^2}
		-
		\frac{1}{4}
		\right],
		\qquad c_1\in[0,1].
		\label{eq:convex_front}
	\end{equation}
	The concave front is
	\begin{equation}
		c_2=g_{\mathrm{conc}}(c_1)
		=
		1-(1-a)c_1-a c_1^2,
		\qquad a=0.5,
		\qquad c_1\in[0,1].
		\label{eq:concave_front}
	\end{equation}
	
	Both fronts satisfy $g(0)=1$ and $g(1)=0$. The concave front contains
	interior non-supported Pareto points with respect to linear scalarization.
	A related preprint studies reachability and front population for different
	scalar objective functions on these fronts~\cite{Frommann2026OScalarObjectives}.
	
	\subsection{Pareto-front-only diagnostic study}
	\label{sec:front_only_study}
	
	In the first study, the fuzzy objectives are evaluated only on the known
	Pareto-front points. The reference point is placed directly on the front,
	
	\begin{equation}
		c_1^*=\lambda,\qquad c_2^*=g(\lambda),
		\qquad \lambda\in\{0.01,0.02,\ldots,0.99\}.
	\end{equation}
	
	This setting is diagnostic rather than directly practical. Since the reference
	point is already Pareto-optimal, exact reference following is not necessarily a
	useful decision objective. The purpose is instead to reveal how memberships and
	rules populate different regions of the front and whether apparent reference
	following is genuine or caused by plateaus.
	
	\subsection{Two-dimensional dominated-reference study}
	\label{sec:2d_study}
	
	To obtain a more practical setting, the same fronts are embedded in a
	two-dimensional design space,
	
	\begin{equation}
		x_1\in[0,1],
		\qquad
		x_2\in[0,1].
	\end{equation}
	
	The first variable controls the position along the Pareto front, while the
	second variable controls the distance into the dominated region. With
	
	\begin{equation}
		\eta=\eta_{\max}x_2,
		\qquad \eta_{\max}=0.5,
	\end{equation}
	
	the criteria are defined as
	
	\begin{align}
		c_1(x_1,x_2)
		&=
		x_1+\eta(1-x_1),\\
		c_2(x_1,x_2)
		&=
		g(x_1)+\eta\bigl(1-g(x_1)\bigr).
	\end{align}
	
\begin{figure}[ht]
		\centering
		\includegraphics[width=\textwidth]{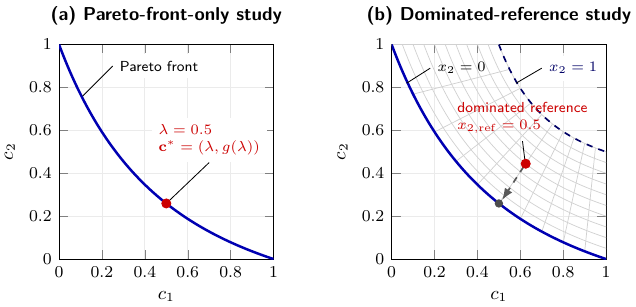}
		\caption{Analytical test setup for the convex front at $\lambda=0.5$.
			The front-only study places the reference on the Pareto front, whereas the
			dominated-reference study embeds the same front in a two-dimensional design
			space with $x_{2,\mathrm{ref}}=0.5$ and Pareto set $x_2=0$.}
		\label{fig:test_setup}
	\end{figure}

	For $x_2=0$, the point lies on the Pareto front. For $x_2>0$, the point is
	dominated by the corresponding point $(x_1,0)$. Thus, the Pareto set is known
	analytically and is given by
	\begin{equation}
		x_2=0.
	\end{equation}
	Figure~\ref{fig:test_setup} illustrates both study setups for the convex front.

	The reference design is chosen as
	
	\begin{equation}
		x_{1,\mathrm{ref}}=\lambda,
		\qquad
		x_{2,\mathrm{ref}}=0.5.
	\end{equation}
	
	This reference design is dominated. The two-dimensional study therefore tests
	whether a fuzzy objective improves a dominated reference design, locks onto it,
	or ignores it in favor of a global Pareto tradeoff.
	
	\section{Evaluation metrics and behavior classes}
	\label{sec:metrics}
	
	The evaluation is designed to identify not only the selected optimum of a fuzzy
	scalar objective, but also the mechanism by which this optimum is selected. 
	This distinction is important because a fuzzy objective may create a localized
	minimum with low tie ambiguity, a broad flat plateau of equally preferred
	points, a fixed compromise independent of the reference, or an endpoint
	attractor. Therefore, the analysis
	records both optimum-dependent quantities and diagnostic quantities such as tie
	fractions, validity fractions, normalized selection entropy, Pareto recovery, and reference
	stagnation.
	
	For each combination of front geometry, membership configuration, rule set, and
	reference value, the fuzzy objective is evaluated on a sampled candidate set. The
	reference parameter is sampled as
	
	\begin{equation}
		\lambda \in \{0.01,0.02,\ldots,0.99\},
	\end{equation}
	
	giving 99 reference cases for each front/configuration/rule-set combination. For
	the front-only study, the candidate set consists of $10001$ uniformly distributed
	points on the analytical Pareto front. For the two-dimensional
	dominated-reference study, a tensor-product grid with $1001$ points in $x_1$ and
	$201$ points in $x_2$ is used.
	
	Let $\mathcal{A}$ denote the sampled candidate set. A candidate is considered
	valid only if the fuzzy objective evaluation returns a defined objective value.
	In particular, candidates for which the denominator in
	Eq.~\eqref{eq:fuzzy_objective} vanishes are excluded from the minimization. The
	valid candidate set is denoted by $\mathcal{A}_{\mathrm{val}}$, and its relative
	size is recorded as
	
	\begin{equation}
		\phi_{\mathrm{val}}
		=
		\frac{|\mathcal{A}_{\mathrm{val}}|}{|\mathcal{A}|}.
		\label{eq:valid_fraction}
	\end{equation}
	
	Cases with no valid candidate point are marked as invalid and are excluded from
	optimum-dependent averages.
	
	For each valid case, the minimum fuzzy objective value is
	
	\begin{equation}
		f_{\min}
		=
		\min_{z\in\mathcal{A}_{\mathrm{val}}}
		f_{\mathrm F}(z),
	\end{equation}
	
	where $z$ denotes either a point on the front or a point in the two-dimensional
	design space. Tied minimizers are detected with the numerical tolerance
	
	\begin{equation}
		\varepsilon_f = 10^{-12}.
	\end{equation}
	
	The set of tied minimizers is therefore
	
	\begin{equation}
		\mathcal{M}
		=
		\left\{
		z\in\mathcal{A}_{\mathrm{val}}
		\;\middle|\;
		|f_{\mathrm F}(z)-f_{\min}|\le \varepsilon_f
		\right\}.
		\label{eq:minimizer_set}
	\end{equation}
	
	The tie count is $|\mathcal{M}|$, and the tie fraction is
	
	\begin{equation}
		\tau
		=
		\frac{|\mathcal{M}|}{|\mathcal{A}_{\mathrm{val}}|}.
		\label{eq:tie_fraction}
	\end{equation}
	
	High tie fractions are interpreted in light of the rule-activation plateau
	mechanism described in Sec.~\ref{sec:rule_activation_plateaus}: if only one
	output class determines the objective over a finite region, the centroid value
	and its scaled objective value become independent of the activation strength.
	All points in that region therefore receive the same scalar objective value.
	
	If several candidates are tied, the reported optimum is chosen as the tied
	candidate closest to the reference point in criterion space,
	
	\begin{equation}
		z_{\mathrm{opt}}
		=
		\arg\min_{z\in\mathcal{M}}
		\left[
		\left(c_1(z)-c_1^*\right)^2
		+
		\left(c_2(z)-c_2^*\right)^2
		\right].
		\label{eq:tie_breaking}
	\end{equation}
	
	This deterministic tie-breaking rule avoids artifacts caused by grid ordering.
	However, it does not remove the ambiguity of the fuzzy objective itself.
	Consequently, the tie fraction is used as a central diagnostic. Low tie
	fractions indicate a localized preference minimum, whereas high tie fractions
	indicate plateau-dominated behavior.
	
	\subsection{Front-only metrics}
	\label{sec:front_only_metrics}
	
	In the front-only study, each candidate point lies directly on the analytical
	Pareto front,
	
	\begin{equation}
		z=(c_1,g(c_1)),
		\qquad c_1\in[0,1].
	\end{equation}
	
	For a given reference value $\lambda$, the reference point is
	
	\begin{equation}
		\mathbf{c}^*
		=
		(c_1^*,c_2^*)
		=
		(\lambda,g(\lambda)).
	\end{equation}
	
	The selected point is summarized by its first criterion value
	$c_{1,\mathrm{opt}}$. The reference error is
	
	\begin{equation}
		e_{\mathrm{ref}}
		=
		|c_{1,\mathrm{opt}}-\lambda|.
		\label{eq:front_reference_error}
	\end{equation}
	
	Both the mean reference error and the root-mean-square reference error are
	computed over all valid reference values. These metrics quantify whether the
	selected point follows the moving reference point along the front.
	
	To detect collapse to a fixed central compromise, a balanced point is computed
	for each front. It is defined as the point at which both normalized criteria are
	equal,
	
	\begin{equation}
		g(c_{1,\mathrm{bal}})=c_{1,\mathrm{bal}}.
		\label{eq:balanced_point_definition}
	\end{equation}
	
	This equation is solved numerically by bisection on $[0,1]$. The balanced-point
	error is
	
	\begin{equation}
		e_{\mathrm{bal}}
		=
		|c_{1,\mathrm{opt}}-c_{1,\mathrm{bal}}|.
		\label{eq:balanced_error}
	\end{equation}
	
	A small mean balanced-point error indicates that a configuration repeatedly
	selects the same central compromise instead of responding to the reference
	parameter.
	
	Endpoint behavior is measured by the indicator
	
	\begin{equation}
		E
		=
		\begin{cases}
			1,
			&
			c_{1,\mathrm{opt}}\le \varepsilon_{\mathrm{end}}
			\;\text{or}\;
			c_{1,\mathrm{opt}}\ge 1-\varepsilon_{\mathrm{end}},
			\\
			0,
			&
			\text{otherwise},
		\end{cases}
		\label{eq:endpoint_indicator_front}
	\end{equation}
	
	with
	
	\begin{equation}
		\varepsilon_{\mathrm{end}}=0.05.
	\end{equation}
	
	The endpoint fraction is the mean value of $E$ over all valid reference values.
	
	The spread of selected front locations is quantified from a histogram of the
	selected values $c_{1,\mathrm{opt}}$. The resulting normalized Shannon entropy
	is referred to as \emph{normalized selection entropy}. It is not used as a
	measure of uncertainty in the fuzzy inference. Instead, it measures how broadly
	the selected front locations are distributed over the front when the reference
	value is varied. The interval $[0,1]$ is divided into
	
	\begin{equation}
		B=10
	\end{equation}
	
	uniform bins. Let $p_b$ be the fraction of selected points in bin $b$. The
	selection entropy is defined as the Shannon entropy~\cite{Shannon1948} of this
	empirical bin distribution,
	\begin{equation}
		H
		=
		-
		\sum_{b=1}^{B}
		p_b \log p_b,
		\label{eq:selection_entropy}
	\end{equation}
	
	where empty bins are omitted from the sum. The normalized selection entropy is
	
	\begin{equation}
		H_{\mathrm{norm}}
		=
		\frac{H}{\log B}.
		\label{eq:normalized_entropy}
	\end{equation}
	
	The logarithm base is immaterial for $H_{\mathrm{norm}}$, provided that the
	same base is used in numerator and denominator.
	Values close to zero indicate that selected points collapse into a small region
	of the front. Values close to one indicate that selected points are distributed
	broadly over the front. The dominant-bin fraction,
	
	\begin{equation}
		p_{\max}
		=
		\max_b p_b,
		\label{eq:dominant_bin_fraction}
	\end{equation}
	
	is used as a complementary concentration measure.
	
	The front-only classes are assigned by the ordered criteria in
	Table~\ref{tab:front_only_class_rules}. The order is part of the definition:
	for example, endpoint collapse is tested before all other classes.
	
	\begin{table}[htbp]
		\centering
		\caption{Front-only behavior classes and assignment criteria. The rules are
			applied in the listed order.}
		\label{tab:front_only_class_rules}
		\begin{tabular}{lll}
			\toprule
			Class & Criterion & Interpretation \\
			\midrule
			EPC & endpoint fraction $>0.90$ & endpoint collapse \\
			BPC & $H_{\mathrm{norm}}<0.05$, $p_{\max}>0.95$,
			$\bar e_{\mathrm{bal}}<0.02$ & balanced-point collapse \\
			PC & $H_{\mathrm{norm}}<0.05$, $p_{\max}>0.95$ & point collapse \\
			ARP & $\bar e_{\mathrm{ref}}<0.02$, $\bar\tau>0.50$ &
			apparent reference-following plateau \\
			GRF & $\bar e_{\mathrm{ref}}<0.02$, $\bar\tau<0.01$,
			$H_{\mathrm{norm}}>0.85$ & genuine reference following \\
			BD & $H_{\mathrm{norm}}>0.80$ & broad front distribution \\
			MIX & otherwise & localized or mixed behavior \\
			\bottomrule
		\end{tabular}
	\end{table}
	
	\subsection{Dominated-reference metrics}
	\label{sec:dominated_reference_metrics}
	
	In the two-dimensional dominated-reference study, the selected point is a design
	point
	
	\begin{equation}
		z_{\mathrm{opt}}
		=
		(x_{1,\mathrm{opt}},x_{2,\mathrm{opt}})
	\end{equation}
	
	with associated criteria
	
	\begin{equation}
		\mathbf{c}_{\mathrm{opt}}
		=
		(c_{1,\mathrm{opt}},c_{2,\mathrm{opt}}).
	\end{equation}
	
	The reference design is
	
	\begin{equation}
		z_{\mathrm{ref}}
		=
		(x_{1,\mathrm{ref}},x_{2,\mathrm{ref}})
		=
		(\lambda,0.5),
	\end{equation}
	
	and its criterion vector is
	
	\begin{equation}
		\mathbf{c}^*
		=
		(c_1^*,c_2^*).
	\end{equation}
	
	Because $x_{2,\mathrm{ref}}=0.5>0$, the reference design is dominated by the
	corresponding design with the same $x_1$ value and $x_2=0$. Therefore, the
	central question in this study is whether the fuzzy scalar objective improves
	the dominated reference, locks onto it, or ignores it in favor of a global
	tradeoff.
	
	Dominance of the reference is tested with the numerical tolerance
	
	\begin{equation}
		\varepsilon_{\mathrm{dom}}=10^{-12}.
	\end{equation}
	
	The selected design is said to dominate the reference if
	
	\begin{equation}
		c_{i,\mathrm{opt}}
		\le
		c_i^*+\varepsilon_{\mathrm{dom}}
		\qquad
		\text{for } i=1,2,
		\label{eq:dominance_no_worse}
	\end{equation}
	
	and if at least one criterion is strictly improved,
	
	\begin{equation}
		c_{j,\mathrm{opt}}
		<
		c_j^*-\varepsilon_{\mathrm{dom}}
		\qquad
		\text{for at least one } j\in\{1,2\}.
		\label{eq:dominance_strict}
	\end{equation}
	
	The dominance rate is the fraction of valid reference cases satisfying
	Eqs.~\eqref{eq:dominance_no_worse} and~\eqref{eq:dominance_strict}.
	
	The total reference violation is
	
	\begin{equation}
		V
		=
		\max(0,c_{1,\mathrm{opt}}-c_1^*)
		+
		\max(0,c_{2,\mathrm{opt}}-c_2^*),
		\label{eq:reference_violation}
	\end{equation}
	
	and the total reference improvement is
	
	\begin{equation}
		I
		=
		\max(0,c_1^*-c_{1,\mathrm{opt}})
		+
		\max(0,c_2^*-c_{2,\mathrm{opt}}).
		\label{eq:reference_improvement}
	\end{equation}
	
	The violation metric measures how much the selected design worsens one or more
	criteria relative to the reference. The improvement metric measures the total
	criterion-wise gain relative to the reference.
	
	Pareto recovery is measured directly in design space because the Pareto set of
	the embedded problem is known analytically. A selected design is counted as
	Pareto recovered if
	
	\begin{equation}
		x_{2,\mathrm{opt}}
		\le
		\varepsilon_{\mathrm P},
		\qquad
		\varepsilon_{\mathrm P}=0.005.
		\label{eq:pareto_recovery}
	\end{equation}
	
	The Pareto-recovery rate is the fraction of valid reference cases satisfying
	Eq.~\eqref{eq:pareto_recovery}.
	
	Reference stagnation is also measured in design space. A selected design is
	counted as remaining at the dominated reference if
	
	\begin{equation}
		|x_{1,\mathrm{opt}}-x_{1,\mathrm{ref}}|
		\le
		\varepsilon_{\mathrm R}
		\quad\text{and}\quad
		|x_{2,\mathrm{opt}}-x_{2,\mathrm{ref}}|
		\le
		\varepsilon_{\mathrm R},
		\label{eq:reference_stagnation}
	\end{equation}
	
	with
	
	\begin{equation}
		\varepsilon_{\mathrm R}=0.005.
	\end{equation}
	
	Endpoint or extreme-tradeoff behavior is detected by
	
	\begin{equation}
		x_{1,\mathrm{opt}}\le 0.05
		\quad\text{or}\quad
		x_{1,\mathrm{opt}}\ge 0.95.
		\label{eq:endpoint_indicator_2d}
	\end{equation}
	
	The tie fraction is again used to separate plateau-driven behavior from
	localized minima. The thresholds are
	
	\begin{equation}
		\tau_{\mathrm{high}}=0.5,
		\qquad
		\tau_{\mathrm{low}}=0.01.
		\label{eq:tie_thresholds}
	\end{equation}
	
	A group is considered high-tie if its mean tie fraction is at least
	$\tau_{\mathrm{high}}$, and low-tie if its mean tie fraction is at most
	$\tau_{\mathrm{low}}$.
	
	The dominated-reference behavior classes are assigned by the following ordered
	rule set. All rates and means are computed over valid reference cases only. Let
	
	\begin{align}
		\bar{\tau} &=
		\text{mean tie fraction},\nonumber \\
		r_{\mathrm P} &=
		\text{Pareto-recovery fraction},\nonumber \\
		r_{\mathrm S} &=
		\text{reference-stagnation fraction},\nonumber \\
		r_{\mathrm D} &=
		\text{dominance rate},\nonumber \\
		r_{\mathrm E} &=
		\text{endpoint fraction},\nonumber \\
		\bar{V} &=
		\text{mean reference violation},\nonumber \\
		\bar{I} &=
		\text{mean reference improvement}.
	\end{align}
	
	The logical conditions \emph{highTie}, \emph{lowTie}, \emph{notPlateau},
	\emph{nearReference}, \emph{dominatesReference}, \emph{paretoRecovered},
	\emph{noViolation}, and \emph{noImprovement} are defined as $\bar\tau\ge0.5$, $\bar\tau\le0.01$, $\bar\tau<0.5$,
	$r_S\ge0.95$, $r_D\ge0.95$, $r_P\ge0.95$, $\bar V\le0.001$,
	and $\bar I\le0.01$, respectively.
	
	\begin{table}[htbp]
		\centering
		\caption{Dominated-reference behavior classes and ordered assignment criteria.}
		\label{tab:dominated_class_rules}
		\begin{tabularx}{\textwidth}{@{}lXl@{}}
			\toprule
			Class & Criterion & Interpretation \\
			\midrule
			RSP &
			nearReference $\wedge$ noImprovement $\wedge$ highTie &
			reference-stagnation plateau \\
			
			RL &
			nearReference $\wedge$ noImprovement &
			reference lock \\
			
			PRI &
			dominatesReference $\wedge$ noViolation $\wedge$ paretoRecovered
			$\wedge$ lowTie &
			Pareto reference improvement \\
			
			WRI &
			dominatesReference $\wedge$ noViolation $\wedge$ notPlateau &
			weak reference improvement \\
			
			GPT &
			paretoRecovered $\wedge$ lowTie &
			global Pareto tradeoff \\
			
			EXT &
			$r_E>0.90$ &
			endpoint or extreme tradeoff \\
			
			MIX &
			otherwise &
			mixed or other behavior \\
			\bottomrule
		\end{tabularx}
	\end{table}
	
	This classification is intentionally ordered. Reference stagnation and reference
	lock are identified before improvement classes. Pareto reference improvement is
	identified before weak reference improvement, and both are identified before
	global Pareto tradeoff. This ensures that configurations are classified
	according to the strongest applicable interpretation of their decision behavior.
	The thresholds are not tuned separately for individual configurations; they are
	fixed once and applied uniformly to all cases.
	
	\section{Results}
	\label{sec:results}
	
	Before discussing individual mechanisms, Table~\ref{tab:behavior_class_counts}
	summarizes the number of behavior-class assignments in both studies. The
	front-only study contains a mixture of collapse, plateau-driven reference
	following, and genuine reference following. In the dominated-reference study,
	global Pareto tradeoffs and reference-stagnation plateaus are the most frequent
	classes.
	
	\begin{table}[htbp]
		\centering
		\caption{Number of observed behavior-class assignments in the two studies. Each
			front contains 24 configuration/rule-set combinations. Classes with zero total assignments are omitted.}
		\label{tab:behavior_class_counts}
		\vspace{0.5em}
		\begin{tabular}{llrrrrrr}
			\toprule
			Study & Front & BPC & ARP & GRF & MIX & EPC & Total \\
			\midrule
			Front-only & convex  & 7 & 7 & 3 & 7 & 0 & 24 \\
			Front-only & concave & 7 & 7 & 9 & 0 & 1 & 24 \\
			\midrule
			\multicolumn{2}{l}{Total} & 14 & 14 & 12 & 7 & 1 & 48 \\
			\bottomrule
		\end{tabular}
		
		\vspace{1em}
		
		\begin{tabular}{llrrrrrr}
			\toprule
			Study & Front & GPT & RSP & RL & WRI & PRI & Total \\
			\midrule
			Dominated-reference & convex  & 15 & 6 & 1 & 1 & 1 & 24 \\
			Dominated-reference & concave & 14 & 6 & 1 & 1 & 2 & 24 \\
			\midrule
			\multicolumn{2}{l}{Total} & 29 & 12 & 2 & 2 & 3 & 48 \\
			\bottomrule
		\end{tabular}
		
		\vspace{0.5em}
		
		{\footnotesize
			BPC: balanced-point collapse; ARP: apparent reference-following plateau;
			GRF: genuine reference following; EPC: endpoint collapse; MIX: localized or
			mixed behavior; GPT: global Pareto tradeoff; RSP: reference-stagnation plateau;
			RL: reference lock; WRI: weak reference improvement; PRI: Pareto reference
			improvement.}
	\end{table}
	
	\subsection{Pareto-front-only diagnostic results}
	\label{sec:front_only_results}
	
\begin{figure}[ht]
		\centering
		\includegraphics[width=0.7\textwidth]{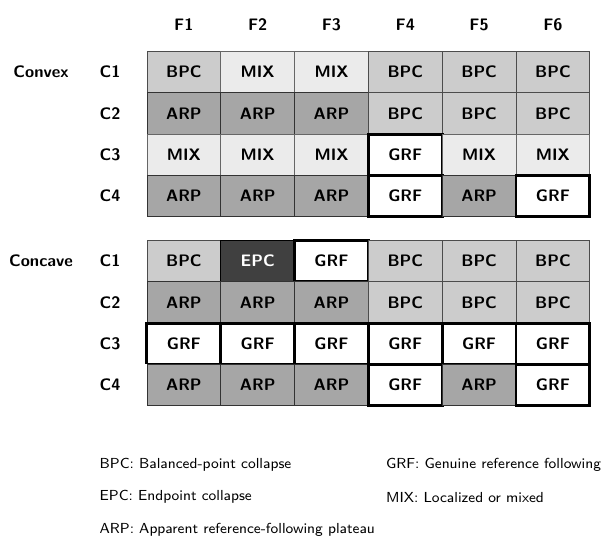}
		\caption{Behavior classification for the Pareto-front-only diagnostic study.
			Each entry aggregates 99 reference values for one front/configuration/rule-set
			combination. Class labels follow the criteria defined in
			Sec.~\ref{sec:front_only_metrics}.}
		\label{fig:front_only_behavior_matrix}
	\end{figure}

	The front-only study provides a diagnostic view of how the fuzzy rule bases
	shape the preference landscape on a known Pareto front. The corresponding
	behavior classification is shown in
	Fig.~\ref{fig:front_only_behavior_matrix}. The metric values reveal that similar
	selected locations can arise from fundamentally different mechanisms.
	
	A first dominant pattern is balanced-point collapse. In the front-only study,
	this behavior is primarily associated with global desirable/undesirable
	membership placement. It occurs for C1 with F1, F4, F5, and F6, and for C2 with
	the multi-consequent rule sets F4--F6. In contrast, the corresponding
	reference-based explicit configuration C4 does not collapse to the balanced
	point; depending on the rule set, it produces either plateau-driven apparent
	reference following or genuine reference following. For example, on the convex
	front, C1 with F1, F4, F5, and F6 selects the balanced point almost exactly: the
	mean balanced-point error is $2.86\times 10^{-5}$, the normalized selection entropy is
	zero, and the dominant-bin fraction is one. The same collapse occurs for
	C2-F4--F6 and is also observed on the concave front, where the balanced-point
	errors are of order $5\times 10^{-5}$. These values show that the selected point
	is essentially independent of the reference parameter. The fuzzy objective
	therefore produces a stable central tradeoff, but does not scan the front with
	the reference point.
	
	A second pattern is apparent reference following caused by flat
	rule-activation plateaus. In the
	explicit-rules-only configurations C2 and C4, the single-consequent rule sets
	F1--F3 often yield mean reference errors equal or close to zero. However, their
	mean tie fractions are equal to one. Thus, all valid sampled points are tied at
	the minimum objective value, and the reported reference-following behavior is a
	consequence of the reference-distance tie-breaking rule. This is visible on both
	fronts. For instance, convex C2-F1, C2-F2, and C2-F3 have zero mean reference
	error but a mean tie fraction of one; the same is observed for the corresponding
	concave cases. These cases demonstrate that small reference error alone is not a
	sufficient indicator of meaningful reference-following behavior.
	
	Genuine reference following requires both small reference error and low tie
	ambiguity. This behavior is obtained most clearly for reference-based
	memberships with competing consequents. For example, C4-F4 and C4-F6 achieve
	zero mean reference error on both fronts, while the normalized selection entropy is close
	to one and the mean tie fraction is approximately $10^{-4}$, corresponding to a
	single selected grid point among the $10001$ sampled front points. In these
	cases, the selected point follows the reference over the front because the fuzzy
	objective forms a localized minimum rather than a broad plateau.
	
	The results also show that the front geometry affects some configurations. In
	particular, C3 produces mostly localized or mixed behavior on the convex front,
	with C3-F4 being the only genuine reference-following case. On the concave
	front, by contrast, all C3 rule sets are classified as genuine reference
	following. This indicates that the same membership/rule configuration can
	interact differently with convex and concave Pareto-front geometries. 
	Table~\ref{tab:front_only_representative_metrics}
	summarizes representative front-only metric values.
	
	\begin{table}[htbp]
		\centering
		\scriptsize
		\caption{Representative front-only metric values. The examples illustrate how
			the behavior classes are supported by the underlying metrics.}
		\label{tab:front_only_representative_metrics}
		\vspace{0.5em}
		\setlength{\tabcolsep}{5.2pt}
		\begin{tabular}{lllrrrrrr}
			\toprule
			Front & Case & Class &
			$\bar e_{\mathrm{ref}}$ &
			$\bar e_{\mathrm{bal}}$ &
			$H_{\mathrm{norm}}$ &
			$p_{\max}$ &
			Endpoint &
			$\bar{\tau}$ \\
			\midrule
			convex  & C1-F1 & BPC & 0.2637 & $2.86{\times}10^{-5}$ & 0.0000 & 1.000 & 0.000 & $1.00{\times}10^{-4}$ \\
			convex  & C2-F1 & ARP & 0.0000 & 0.2637 & 0.9998 & 0.101 & 0.101 & 1.000 \\
			convex  & C4-F6 & GRF & 0.0000 & 0.2637 & 0.9998 & 0.101 & 0.101 & $1.00{\times}10^{-4}$ \\
			convex  & C3-F3 & MIX & 0.1619 & 0.1018 & 0.5920 & 0.313 & 0.000 & $1.02{\times}10^{-4}$ \\
			concave & C1-F2 & EPC & 0.2568 & 0.5081 & 0.2973 & 0.566 & 1.000 & $2.00{\times}10^{-4}$ \\
			\bottomrule
		\end{tabular}
		
		\vspace{0.5em}
		\footnotesize
		$H_{\mathrm{norm}}$: normalized selection entropy;
		$p_{\max}$: dominant-bin fraction;
		$\bar{\tau}$: mean tie fraction.
	\end{table}
	
	Overall, the front-only study supports three conclusions. First, global
	memberships with implicit baseline rules tend to produce stable central
	tradeoffs rather than reference-following behavior. Second, explicit-rules-only
	single-consequent rule bases can produce apparent reference following through
	flat plateaus. Third, low tie fractions are essential for distinguishing genuine
	reference-following behavior from tie-breaking artifacts.
	
	\subsection{Dominated-reference results}
	\label{sec:dominated_reference_results}
	
	The two-dimensional dominated-reference study changes the interpretation of the
	reference point. The reference design is now intentionally dominated, with
	$x_{2,\mathrm{ref}}=0.5$, while the Pareto set is located at $x_2=0$. A useful
	reference-improvement objective should therefore move away from the reference
	toward smaller $x_2$, ideally reaching the Pareto set without worsening either
	criterion relative to the reference.
	
\begin{figure}[ht]
		\centering
		\includegraphics[width=0.7\textwidth]{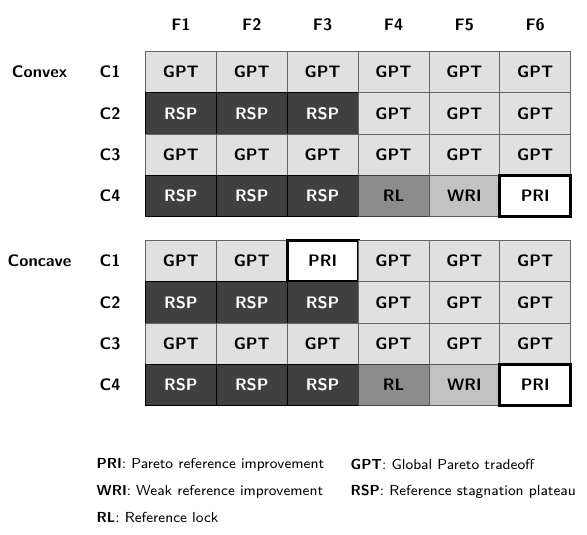}
		\caption{Behavior classification for the two-dimensional dominated-reference
			study. Each entry aggregates 99 dominated reference designs for one
			front/configuration/rule-set combination. Class labels follow the ordered
			criteria defined in Sec.~\ref{sec:dominated_reference_metrics}.}
		\label{fig:2d_behavior_matrix}
	\end{figure}

	The behavior classification is shown in
	Fig.~\ref{fig:2d_behavior_matrix}. The matrix reveals a highly structured
	pattern across configurations and rule sets.
	
	Configurations with global memberships generally recover the Pareto set but do
	not necessarily dominate the reference for every reference value. These cases
	are classified as global Pareto tradeoffs. For example, convex C1-F1 has a
	Pareto-recovery fraction of one and a mean selected value
	$x_{2,\mathrm{opt}}=0$, but it dominates the reference in only $47.5\%$ of the
	reference cases and has a nonzero mean reference violation of $0.0342$. Similar
	patterns occur for many C1, C2-F4--F6, and C3 cases. These objectives therefore
	produce reliable Pareto solutions with low tie fractions, but their selected
	tradeoff is governed by a global compromise mechanism rather than by systematic
	improvement of each reference design.
	
	The explicit-rules-only single-consequent rule sets F1--F3 behave very
	differently. In C2 and C4 they produce reference stagnation plateaus on both
	fronts. The selected designs remain at the dominated reference, with
	$x_{2,\mathrm{opt}}=0.5$, reference-stagnation fraction equal to one, Pareto
	recovery equal to zero, and mean tie fraction equal to one. Thus, the plateau
	mechanism identified in the front-only study becomes practically harmful in the
	dominated-reference setting: instead of improving the dominated reference, the
	objective remains at it.
	
\begin{figure}[ht]
		\centering
		\includegraphics[width=\textwidth]{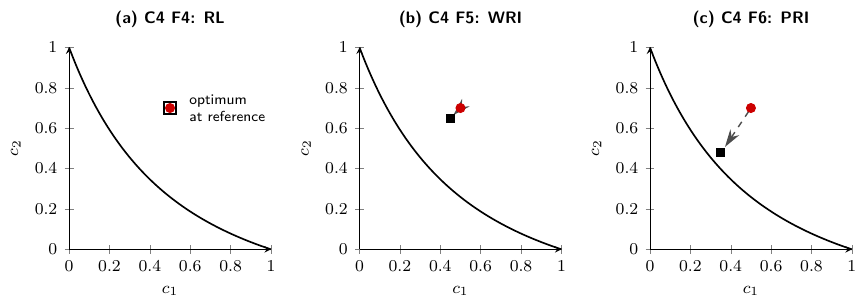}
		\caption{Representative C4 cases in the dominated-reference study. F4 locks
			onto the dominated reference, F5 produces only weak improvement, and F6
			improves the reference while recovering the Pareto set.}
		\label{fig:c4_representative_cases}
	\end{figure}

	The C4 cases are particularly informative because they use reference-based
	memberships with explicit-rules-only inference. The different mechanisms are
	illustrated in Fig.~\ref{fig:c4_representative_cases} and quantified in
	Table~\ref{tab:c4_metric_summary}.
	
	\begin{table}[htbp]
		\centering
		\scriptsize
		\caption{Quantitative summary of the C4 dominated-reference cases. C4 uses
			reference-based memberships with explicit-rules-only inference. The metrics
			show the transition from reference stagnation plateau (RSP) and reference
			lock (RL) to weak reference improvement (WRI) and Pareto reference
			improvement (PRI).}
		\label{tab:c4_metric_summary}
		\vspace{0.5em}
		\setlength{\tabcolsep}{3.3pt}
		\begin{tabular}{llcccccccc}
			\toprule
			Front & Rule & Class & $r_D$ & $\bar V$ & $\bar I$ &
			$\bar{x}_{2}$ & $r_P$ & $r_S$ & $\bar{\tau}$ \\
			\midrule
			convex  & F1--F3 & RSP & 0.000 & $\le 7.5{\times}10^{-4}$ & $\le 7.5{\times}10^{-4}$ & 0.500 & 0.000 & 1.000 & 1.000 \\
			convex  & F4     & RL  & 0.000 & 0.000 & 0.000 & 0.500 & 0.000 & 1.000 & $9.64{\times}10^{-2}$ \\
			convex  & F5     & WRI & 1.000 & 0.000 & $3.04{\times}10^{-3}$ & 0.495 & 0.000 & 0.000 & $9.64{\times}10^{-2}$ \\
			convex  & F6     & PRI & 1.000 & 0.000 & 0.356 & 0.000 & 1.000 & 0.000 & $4.97{\times}10^{-6}$ \\
			\midrule
			concave & F1--F3 & RSP & 0.000 & $\le 7.5{\times}10^{-4}$ & $\le 7.5{\times}10^{-4}$ & 0.500 & 0.000 & 1.000 & 1.000 \\
			concave & F4     & RL  & 0.000 & 0.000 & 0.000 & 0.500 & 0.000 & 1.000 & $6.02{\times}10^{-2}$ \\
			concave & F5     & WRI & 1.000 & 0.000 & $2.05{\times}10^{-3}$ & 0.495 & 0.000 & 0.000 & $6.02{\times}10^{-2}$ \\
			concave & F6     & PRI & 1.000 & 0.000 & 0.205 & 0.000 & 1.000 & 0.000 & $4.97{\times}10^{-6}$ \\
			\bottomrule
		\end{tabular}
		
		\vspace{0.5em}
		\footnotesize
		$r_D$: dominance rate; $\bar V$: mean reference violation;
		$\bar I$: mean reference improvement; $\bar{x}_{2}$: mean selected
		$x_2$ coordinate; $r_P$: Pareto-recovery fraction; $r_S$:
		reference-stagnation fraction; $\bar{\tau}$: mean tie fraction.
	\end{table}
	
	Rule set F4 in C4 produces reference lock rather than plateau-driven
	stagnation. On both fronts, the selected design remains at the reference
	design, with reference-stagnation fraction equal to one and Pareto recovery
	equal to zero. However, the mean tie fraction is lower than in the stagnation
	plateau cases: approximately $0.096$ on the convex front and $0.060$ on the
	concave front. This indicates that the lock is not caused by a full flat
	plateau. Instead, it follows from the rule semantics of F4, which maps the
	simultaneous tolerable reference neighborhood directly to the desirable output
	class. The reference region itself therefore becomes the preferred outcome.
	
	Rule set F5 in C4 produces weak reference improvement. It dominates the
	reference in all tested reference cases and has zero mean reference violation on
	both fronts. However, the improvement is very small and the Pareto set is not
	recovered. On the convex front, the mean improvement is only $0.0030$ and the
	mean selected value is $x_{2,\mathrm{opt}}=0.495$. On the concave front, the
	mean improvement is $0.0020$, again with
	$x_{2,\mathrm{opt}}=0.495$. Hence, F5 moves the design only slightly away from
	the dominated reference and does not provide robust Pareto recovery.
	
	The strongest and most consistent reference-improvement behavior is obtained
	with C4-F6. On both the convex and concave fronts, C4-F6 is classified as
	Pareto reference improvement. It dominates the reference in all reference cases,
	has zero mean reference violation, recovers the Pareto set in all cases, and has
	a very low mean tie fraction of approximately $5\times 10^{-6}$. The mean
	selected value is $x_{2,\mathrm{opt}}=0$ on both fronts. The mean reference
	improvement is $0.356$ on the convex front and $0.205$ on the concave front.
	This confirms that the three-class rule structure of F6 is essential in the
	dominated-reference setting: better-than-reference regions can become desirable,
	the reference neighborhood remains only tolerable, and worse regions are
	penalized as undesirable.
	
	One additional front-dependent exception is observed for C1-F3 on the concave
	front. This configuration is also classified as Pareto reference improvement,
	with dominance rate one, zero mean violation, Pareto recovery equal to one, and
	low tie fraction. However, the same configuration is classified as a global
	Pareto tradeoff on the convex front. Thus, although C1-F3 can perform well for a
	specific front geometry, it does not provide the same geometry-independent
	reference-improvement behavior as C4-F6.
	
	In summary, the dominated-reference study separates three practically different
	uses of fuzzy scalar objectives. Global configurations with competing
	consequents provide robust Pareto tradeoffs, but do not necessarily improve each
	reference design. Reference-attraction rules can intentionally lock onto a
	reference region, which is useful only if the reference is meant to be a target.
	For dominated-reference improvement, C4-F6 is the most reliable configuration
	in the present tests because its three-class rule semantics distinguish
	better-than-reference, reference-neighborhood, and worse-than-reference regions.
	
	\section{Discussion}
	\label{sec:discussion}
	
	The results show that the role of the reference point must be specified before
	a fuzzy scalar objective can be designed reliably. Reference following is not a
	universal goal. If the reference point represents a target region on the desired
	Pareto front, then reference-locking or reference-attracting behavior may be
	intended. If the reference point represents a dominated baseline design,
	however, the objective should improve the reference rather than lock onto it.
	
	This distinction explains the difference between F4 and F6. F4 maps reference 
	neighborhoods directly to the desirable class. It is therefore
	appropriate when the reference region represents the intended target. In the
	dominated-reference study, however, this semantics causes F4 to lock onto the
	dominated reference instead of improving it. F6 preserves the
	three-class structure. Reference-neighborhood points are only tolerable, while
	better-than-reference points can become desirable. This makes F6 more suitable
	for improving dominated reference designs.
	
	The study also shows that explicit-rules-only inference is not automatically
	reliable. Single-consequent rule bases often produce flat rule-activation plateaus. 
	In such cases, the selected point may appear to follow the reference,
	but the choice is controlled by tie-breaking rather than by a meaningful
	preference gradient. The tie fraction is therefore an essential diagnostic
	quantity. High tie fractions indicate that the rule base is not sufficiently
	discriminating.
	
	Global memberships remain useful when the decision intent is not reference
	improvement but robust compromise selection. In these cases, global
	desirable/undesirable memberships with competing rule consequents produce
	Pareto tradeoffs with low tie fractions. Such configurations are appropriate 
	when the user has no reliable reference design or wants to avoid severely 
	poor individual criterion values.
	
	These findings have a direct implication for fuzzy decision models used as
	optimization objectives. A linguistically plausible rule base is not sufficient;
	the induced scalar landscape must also be checked. In particular, rule bases
	with only one active consequent class may express a meaningful qualitative
	statement, but after normalized aggregation they can become non-discriminating
	over large regions of the candidate set. Therefore, fuzzy decision models used
	inside optimization loops should be diagnosed not only by inspecting membership
	functions and rules, but also by evaluating tie fractions, reference violations,
	and Pareto-recovery behavior.
	
	\subsection{Limitations}
	\label{sec:limitations}
	
	The study is intentionally restricted to controlled bi-criteria problems. This
	allows exact visualization of the Pareto front, dominated reference designs,
	plateaus, and tie-breaking behavior. Higher-dimensional problems may introduce
	additional rule-interaction effects, so the results should be interpreted as
	diagnostic evidence for important mechanisms rather than as a complete taxonomy
	of fuzzy scalarization behavior.
	
	Parabolic memberships are used throughout to isolate the effects of membership
	placement and rule consequents; other membership families may lead to additional
	behavior. Finally, exhaustive sampling removes optimizer-dependent effects and
	exposes the objective landscape, but practical continuous optimization may
	require additional robustness checks.
	
	\subsection{Practical decision guidance}
	\label{sec:guidance}
	
	The results can be summarized as three practical decision questions.
	
	First, if no reliable reference design is available, global desirable and
	undesirable memberships should be used with competing desirable and
	undesirable consequents. Rule sets such as F5 or F6 then provide a robust global
	Pareto tradeoff. This is suitable for conservative compromise selection, but it
	does not guarantee improvement of an arbitrary reference design.
	
	Second, if the reference point is intended as a target region, reference-based
	memberships and attraction rules can be appropriate. In this case, locking onto
	the reference may be desired. Rule set F4 represents this behavior by mapping
	the reference-neighborhood premise directly to the desirable output class.
	
	Third, if the reference point is a dominated baseline or aspiration design, it
	should not be mapped directly to a desirable outcome. The recommended template 
	is a three-class reference-based rule base such as F6
	in combination with reference-based membership functions: desirable--desirable
	combinations map to desirable, tolerable--tolerable combinations map to
	tolerable, and undesirable deviations map to undesirable. In the present tests, 
	this is the only configuration that
	consistently improves the dominated reference, recovers the Pareto set, and
	avoids plateau-driven selection on both test fronts.
	
	\section{Conclusions}
	\label{sec:conclusions}
	
	This paper investigated how membership placement, implicit rules, and rule
	consequents affect fuzzy scalar objective functions for multi-criteria decision
	making. A Pareto-front-only ablation study first isolated rule-induced
	front-selection mechanisms. It showed that reference following can be genuine
	or apparent, and that apparent reference following may be caused by flat
	rule-activation plateaus.
	
	A second, more practical, two-dimensional dominated-reference study showed how
	these mechanisms affect optimization from a dominated baseline design. The 
	results demonstrate that the appropriate fuzzy configuration depends on the
	role assigned to the reference point. Global memberships with competing consequents yield robust
	global Pareto tradeoffs. Reference-locking rules may be appropriate when the
	reference is a target region. When the reference is a dominated baseline to be
	improved, however, the three-class F6 rule set with reference-based memberships
	provides the most reliable behavior in the present tests: it dominates the
	reference, avoids reference violation, recovers the Pareto set, and has a low tie
	fraction.
	
	The practical conclusion is that fuzzy scalar objectives should be designed from
	the intended decision semantics rather than from membership functions alone. In
	particular, users should distinguish reference targets from dominated reference
	baselines and should diagnose plateau behavior using tie fractions. The present
	results are limited to controlled bi-criteria test cases; whether the same
	patterns generalize to higher-dimensional problems remains to be investigated.

	\section*{Data and code availability}
	
	The source code, generated numerical data, classification scripts, and the
	vendored FuzzyGoal 1.2.0 dependency snapshot are archived on Zenodo at
	\url{https://doi.org/10.5281/zenodo.21310754}.
	
	\section*{Declarations}
	
	No specific funding was received for conducting this study. The author declares
	no relevant competing interests. Generative AI assistance was used for language
	editing, documentation structuring, drafting support, and consistency checks;
	the author reviewed and validated the final manuscript and remains fully
	responsible for its content.
	
	\printbibliography
	
\end{document}